\documentstyle[12pt]{article}

\newcommand{\sect}[1]
{\setcounter{equation}{0}\section{#1}}

\def\noi{\noindent}
\def\ar{\begin{array}{rcl}}
\def\an{\end{array}}

\newcommand{\eq}{\begin{equation}}
\newcommand{\eqa}{\begin{eqnarray}}
\newcommand{\en}{\end{equation}}
\newcommand{\ena}{\end{eqnarray}}
\def\1{{\bf 1}}

\def\ot{\otimes}
\def\id{\mbox{id}}
\def\ie{\mbox{\it i.e.\/ }}
\def\eg{\mbox{\it e.g.\/ }}

\def\g{\mbox{\bf g\,}}
\def\uqg{\mbox{$U_{q}{\/\mbox{\bf g}}$ }}
\def\fun{\mbox{Fun$(G_{q})$\,}}

\def\R{\mbox{$\hat R$\,}}
\def\RH{\mbox{$\hat {\sf R}$\,}}

\def\A{\mbox{${\cal A}_
{\pm,{\tiny G},\phi}$}}
\def\Aq{\mbox{${\cal A}^q_
{\pm,{\footnotesize G\,},\phi}$}}

\def\P{\mbox{$\cal P$}}
\def\PH{\mbox{$\sf P$}}

\def\z{\hspace*{9mm}}
\def\x{\hspace{3mm}}

\newcommand{\rn}{{\bf R}}

\newtheorem{prop}{Proposition}
\newtheorem{lemma}{Lemma}

\begin{document}
\begin{titlepage}
\begin{center}
~

\vskip.6in

{\Large \bf Braided Chains of q-Deformed\\
 Heisenberg Algebrae}

\vskip.4in

Gaetano Fiore*

\vskip.25in

{\em Sektion Physik der Ludwig-Maximilians-Universit\"at
M\"unchen\\
Theoretische Physik --- Lehrstuhl Professor Wess\\
Theresienstra\ss e 37, 80333 M\"unchen\\
Federal Republic of Germany}

\end{center}
\vskip1.6in

\begin{abstract}

Given $M$ copies of a $q$-deformed
Weyl or Clifford algebra in the defining representation of a
quantum group $G_q$, we determine a prescription
to embed  them into a unique, inclusive $G_q$-covariant algebra. 
The different copies are ``coupled" to each other and are
naturally ordered into a ``chain". 
In the case $G_q=SL_q(N)$ a modified
prescription yields an inclusive algebra which is even 
explicitly $SL_q(M)\times SL_q(N)$-covariant,
where $SL_q(M)$ is a
symmetry relating the different copies.
By the introduction of these inclusive algebrae 
we significantly enlarge the class of 
$G_q$-covariant deformed Weyl/Clifford algebrae available
for physical applications.

\end{abstract}

\vfill

\noi \hrule
\vskip.2cm
\noi{\footnotesize *EU-fellow, TMR grant ERBFMBICT960921.
\qquad {\it e-mail: }Gaetano.Fiore@physik.uni-muenchen.de}
\end{titlepage}
\newpage
\setcounter{page}{1}

\renewcommand{\theequation}{\thesection.\arabic{equation}}
\sect{Introduction}
\label{intro}

Weyl and Clifford algebrae (respectively denoted by 
${\cal A}_+,{\cal A}_-$ 
in the sequel, and collectively as ``Heisenberg
algebrae'') are at the hearth of quantum physics.
The most useful Heisenberg algebrae are
the ones endowed with definite
transformation properties under the action
of some symmmetry Lie group $G$ (or Lie algebra \g).

The idea that quantum groups \cite{dr2}
could generalize Lie groups in describing symmetries of
quantum physical systems has attracted much interest
in the past decade. Mathematically speaking,
a quantum group can be described as a deformation
\fun of the algebra $Fun(G)$ of regular functions on $G$ or,
in the dual picture, as a deformation \uqg of
the universal enveloping algebra $U\g$,
within the category of
(quasitriangular) Hopf algebrae; here $q=e^h$, and
$h$ is the deformation parameter.  These 
$q$-deformations induce 
matched $q$-deformations of all 
$\fun$-comodule algebrae [\ie
of the algebrae whose generators satisfy 
commutation relations
that are preserved by the $\fun$-coaction], in 
particular of
$G$-covariant Heisenberg algebrae.
$q$-Deformed Heisenberg algebrae 
corresponding to a 
simple Lie algebra $\g$ in the classical 
series $A_n,B_n,D_n$
were introduced in Ref. \cite{puwo,wezu,pusz,cawa} in the
restricted case that the generators $A^+_i,A^i$ belong
respectively to the defining corepresentation
$\phi_d$ of $\fun$ and to its contragradient 
$\phi_d^\vee$. 

In general,
we shall denote by $\A$ the Weyl/Clifford algebra with
generators $a^i,a^+_i$ belonging respectively to
some corepresentation $\phi$ of $G$
and to its contragradient $\phi^\vee$ and fulfilling the
canonical (anti)commutation relations 
\eqa
a^+_i \,a^+_j \mp a^+_j\, a^+_i &=& 0\label{ccr1} \\
a^i \, a^j \mp a^j \, a^i &= &0\label{ccr2} \\
a^i\, a^+_j \:-\:\delta^i_j\1 \mp a^+_j\, a^i &= &0. \label{ccr3}
\ena
The purpose of this work is to 
find out if there exists 
some $G_q$-covariant deformation of $\A$
(which we will denote by $\Aq$)
having the same
Poincar\'e series as $\A$. We shall denote the
generators of $\Aq$ by
$A^i,A^+_i$.   

As a preliminary result we show (Sect. \ref{qcrsec}) that, beside 
${\cal A}^q_{\pm,SL(N)\,,\phi_d}$\footnote{$SL(N)$ can be easily
promoted also to a $GL(N)$},
${\cal A}^q_{+,SO(N)\,,\phi_d}$
\cite{puwo,wezu,pusz,cawa}, also 
${\cal A}^q_{-,Sp(n)\,,\phi_d}$
can be defined. The first major result is
however that one can embed $M$ identical copies
of  ${\cal A}^q_{+,G\,,\phi_d}$ (resp.
${\cal A}^q_{-,G\,,\phi_d}$) into a unique,
well-defined algebra
${\cal A}^q_{+,G\,,\phi_M}$ (resp.
${\cal A}^q_{-,G\,,\phi_M}$),  or more generally
$M'<M$ copies
of  ${\cal A}^q_{+,G\,,\phi_d}$ and $(M\!-\!M')$
copies of  ${\cal A}^q_{-,G\,,\phi_d}$ into 
a unique, well-defined deformed superalgebra
${\cal A}^q_{G\,,\phi_M}$; $\phi_M$ 
denotes here the
direct sum of $M$ copies of $\phi_d$. Due to the 
rules of braiding \cite{majid},
the different copies do not commute with each other; 
consistent
commutation relations between the latter require the 
introduction
of an ordering: we call the orderd sequence a
``braided chain".

The use of the symbols $a^i,a^+_i,A^i,A^+_i$ etc. does
not necessarily mean  that we are dealing with 
creators \& annihilators; the latter fact
is rather determined by
the choice of the $*$-structure, if any. In section \ref{star}
we consider the natural $*$-structures
giving the generators the meaning of creation \& 
annihilation  operators, or \eg of
coordinates and derivatives.

The second major result (Sect. \ref{glcov})
is that if $G_q=SL_q(N)$ one can modify the
$A$-$A^+$ commutation relations of
${\cal A}^q_{\pm,SL(N)\,,\phi_M}$ in such a way that the generators
become explicitly $GL_q(M)\times SL_q(N)$-covariant\footnote{The 
result regarding ${\cal A}^q_{\pm,SL(N)\,,\phi_M}$ was essentially 
already found in \cite{quesne}, whose author we thank for drawing 
our attention to this point.}. The additional symmetry
$GL_q(M)$ transforms the different copies into each other,
as in the classical case.

The physical relevance of the case
that $\phi$ is a direct sum of many copies of $\phi_d$'s 
is easily understood once one notes that
the different copies could correspond to different
particles, crystal sites or space(time)-points, 
respectively in quantum mechanics, condensed
matter physics or quantum field theory.
The coupling (\ie non-commutativity) between the different
copies can be interpreted as a naturally built-in form of 
interaction between them.
In the particular case that $\Aq$ (with $q\in \rn$)
is a $q$-deformation of 
the $*$-algebra $\A$ with 
$(a^i)^{\dagger}=a^+_i$, then the 
physical
interpretation of $A^i,A^+_j$ as annihilators and
creators does not necessarily 
requires the introduction
of particles with {\it exotic} statistics. Indeed, it is possible
to adopt ordinary boson/fermion statistics 
\cite{fiogrou,fiocmp}, whereby
$A^+_i,A^i$ are to be interpreted as ``composite operators''
creating and destroying some sort of ``dressed states'' 
of bosons/fermions.  

\sect{Preliminaries}
\label{preli}

For a simple Lie group $G$ the algebra \fun \cite{frt} is generated by
$N^2$ objects $T^i_j$, $i,j=1,...,N$, fulfilling the 
commutation relations
\eq
\R^{ij}_{hk}T^h_lT^k_m= T^i_hT^j_k\R^{hk}_{lm}.
\label{rel1}
\en
$N$ is the dimension of the defining representation of $G$, 
$\R$ the corresponding `braid matrix' \cite{frt}, \ie a numerical matrix
fulfilling the `braid equation' 
\eq
\hat R_{12}\,\hat R_{23}\,\hat R_{12} =\hat R_{23}\,\hat R_{12}\,
\hat R_{23}                                              
\label{braid1}.
\en
Here we have used the conventional tensor notation
$(M_{12})^{ijk}_{lmn}=M^{ij}_{lm}\delta^k_n$, {\it etc}.
Because of eq.'s (\ref{braid1}), (\ref{rel1}) \fun is also a bialgebra
with
coproduct and counit respectively given by
$\Delta(T^i_j)=T^i_h\otimes T^h_j$ and 
$\varepsilon(T^i_j)=\delta^i_j$. 

A (right)  comodule algebra 
of \fun is an algebra ${\cal M}$ equipped with a `corepresentation'
$\phi$, \ie with an algebra
homomorphism $\phi:{\cal M}\rightarrow {\cal M}\otimes\fun $
such that $(\id \otimes\Delta)\circ \phi=(\phi\otimes\id)\circ \phi$.
For any polynomial function $f(t)$ in one variable,
the algebra ${\cal M}$  generated by $N$ objects
$A^+_i$  fulfilling the quadratic
relations
\eq
[f(\R)]^{ij}_{hk}A^+_iA^+_j=0
\label{zaza1}
\en
and equipped with the algebra homomorphism 
$\phi_d(A^+_i):=A^+_j\otimes T^j_i$
is a comodule algebra \cite{frt}.

By adding to the quadratic conditions (\ref{rel1}) some suitable
inhomogeneous condition \cite{frt}, 
\fun can be endowed also with an antipode
$S$ and therefore becomes a Hopf algebra\footnote{
In the case $G_q=SL_q(N)$ this condition reads $\mbox{det}_qT=1$,
where $\mbox{det}_qT$ is the $q$-deformed determinant of $T$.
One can also define a Hopf algebra $GL_q(N)$ by using
the same $\R$-matrix, introducing a new generator $t$ that
is central and group-like, together with
its inverse $t^{-1}$, and then imposing the weaker condition
$\mbox{det}_qT=t$.}. Then the algebra
${\cal M'}$  generated by $N$ objects $A^i$  fulfilling the quadratic
relations
\eq
[f(\R)]_{ij}^{hk}A^jA^i=0
\label{zaza2}
\en
and equipped with the algebra
homomorphism $\phi_d^\vee(A^i):=A^j\otimes ST_j^i$  is a comodule
algebra with inverse transformation properties of ${\cal M}$; therefore 
the corepresentation $\phi_d^\vee$ can be called the contragradient of
$\phi_d$.

To go on, we need to recall some specific information regarding
each quantum group $G_q$.
The braid matrix $\hat R$ of the quantum group $\fun$
is a $N^2\times N^2$ matrix that admits
the following projector decomposition
\cite{frt}
\eq
\begin{array} {lclrcl}
\hat R & = & q{\cal P}^S - q^{-1}{\cal P}^A \qquad & 
\mbox{if }G
& = & SL(N)
    \label{decosl} \cr
\hat R & = & q{\cal P}^s - q^{-1}{\cal P}^a + 
q^{1-N}{\cal P}^t   
\qquad & \mbox{if }G &=& SO(N) \label{decoso} \cr
\hat R & = &
q{\cal P}^{s'} - q^{-1}{\cal P}^{a'} - q^{1-N}{\cal P}^{t'} 
\qquad & \mbox{if }G &=& Sp(n), ~~N=2n, \label{decosp}
\end{array}
\label{decompo}
\en
with 
\eq
{\cal P}^{\mu}{\cal P}^{\nu} = \delta^{\mu\nu}, \qquad
\sum_{\mu}{\cal P}^{\mu} = 1.
\label{orto}
\en
$\P^A,\P^S$ are $SL_q(N)$-covariant 
$q$-deformations of the antisymmetric and  symmetric
projectors respectively; they have dimensions
$\frac{N(N-1)}2$ and $\frac{N(N+1)}2$ respectively. 
${\cal P}^a,{\cal P}^t,{\cal P}^s$ are 
$SO_q(N)$-covariant
$q$-deformations of the antisymmetric, trace, 
and symmetric trace-free projectors respectively; they have dimensions
$\frac{N(N-1)}2,1$ and $\frac{N(N+1)}2-1$ respectively. 
${\cal P}^{s'}, {\cal P}^{t'},{\cal P}^{a'}$ are 
$Sp_q(n)$-covariant ($N=2n$)
$q$-deformations respectively of the symmetric,
symplectic, antisymmetric symplectic-free projectors; 
they have dimensions $\frac{N(N+1)}2$,1 and
$\frac{N(N-1)}2-1$ respectively. Setting 
\eq
\begin{array}{lclrcl}
\P^+ &=& \P^S \qquad\qquad & \mbox{if }G &=& SL(N) \cr
\P^+ &=& \P^s+\P^t \qquad\qquad & \mbox{if }G 
&=& SO(N) \cr
\P^+ &=& \P^{s'} \qquad\qquad & \mbox{if }G &=& Sp(n) \cr
\P^- &=& \P^A \qquad\qquad & \mbox{if }G &=& SL(N) \cr
\P^- &=& \P^a \qquad\qquad & \mbox{if } G &=& SO(N) \cr
\P^- &=& \P^{a'}+\P^{t'} \qquad\qquad &
\mbox{if } G &=& Sp(n) 
\end{array}
\label{defproj}
\en
we obtain \fun-covariant deformations $\P^+,\P^-$ of the
$\frac N2(N\!+\!1)$-dim symmetric 
and $\frac N2(N\!-\!1)$-dim 
antisymmetric projector respectively. 

In the sequel we shall need also the explicit expression
for the $\R$ matrix of $SL_q(N)$ and for its inverse:
\eqa
\R & = & q\sum\limits_{i=1}^N e^i_i\otimes e^i_i
+\sum\limits_{i\neq j} e^j_i\otimes e^i_j+(q-q^{-1})
\sum\limits_{i<j} e^i_i\otimes e^j_j.
\label{explicit}\\
\R^{-1} & = & q^{-1}\sum\limits_{i=1}^N e^i_i\otimes e^i_i
+\sum\limits_{i\neq j} e^j_i\otimes e^i_j+(q^{-1}-q)
\sum\limits_{i>j} e^i_i\otimes e^j_j.
\label{explicit'}
\ena
Here we have used the conventional tensor notation and
denoted by $e^i_j$ the $N\times N$ matrix with 
$(e^i_j)^h_k=\delta^{ih}\delta_{jk}$.

By repeated application of the equations (\ref{braid1}),
(\ref{rel1}) we find
\eq
\begin{array}{rcl}
f(\R)^{ij}_{hk}T^h_lT^k_m &=& T^i_hT^j_kf(\R)^{hk}_{lm},\cr
f(\hat R_{12})\,\hat R_{23}\,\hat R_{12}
&=&\hat R_{23}\,\hat R_{12}\,f(\hat R_{23})                     
\end{array}
\label{braid2}
\en
for any polynomial function $f(t)$ in one variable, in particular for
those $f$'s yielding $f(\hat R)={\cal P}^{\mu}$ or
$f(\hat R)=\hat R^{-1}$.  The Equations~(\ref{braid1}), (\ref{rel1})
and~(\ref{braid2}) are evidently satisfied also after the replacement
$\hat R\rightarrow\hat R^{-1}$.
   
\medskip

If in relations (\ref{zaza1}), (\ref{zaza2}) one chooses
$f(\R)=\P^{\mp}$, then these equations become
the $\fun$-covariant
deformations of the (anti)commutation relations (\ref{ccr1}),
(\ref{ccr2}): 
\eqa
\P^{\mp~ij}_{~~hk}A^+_i\, A^+_j & = & 0, \label{qcr1}\\
\P^{\mp~ij}_{~~hk}A^k\, A^h & = & 0.  \label{qcr2}
\ena
Relations (\ref{ccr1}), (\ref{ccr2}), (\ref{qcr1}), (\ref{qcr2})
amount each to $\frac{N(N-1)}{2}$ or
to $\frac{N(N+1)}{2}$ independent relations, respectively if the upper
or
the lower sign is considered.
The algebrae ${\cal M}$, ${\cal M}'$ defined resp.
by (\ref{qcr1}), (\ref{qcr2}) have
\cite{frt,fiodet} 
 the same Poincar\'e series as the algebrae defined by resp.
by (\ref{ccr1}), (\ref{ccr2}).

To obtain \fun-covariant deformations 
${\cal A}^q_{\pm,G,\phi_d}$ of the classical
Heisenberg algebrae described in Section \ref{intro} one still has to
deform relations (\ref{ccr3}). For $G_q=SL_q(N),SO_q(N)$
this was done in Ref. \cite{puwo,wezu,pusz,cawa}. 
The natural ansatz is to look for
quadratic cross commutation relations, in the form
\eq
A^i\, A^+_j=\delta^i_j{\bf 1} \:\pm\: S^{ih}_{jk}
A^+_h\, A^k.
\label{ansatz}
\en
The inhomogeneous term is fixed by the requirement that $\{A^i\}$
is the basis dual to $\{A^+_i\}$. The numerical
matrix $S$ has to be determined imposing \fun-covariance and that
that ${\cal A}^q_{\pm,G,\phi_d}$ itself has
the same Poincar\'e series as its classical counterpart 
${\cal A}_{\pm,G,\phi_d}$.
It will be convenient to use the following general 

\begin{lemma}
Let $\R=\sum_{\mu}c_{\mu}\P^{\mu}$ 
be the projector decomposition of the braid matrix $\R$, 
and let $\P^+:=\sum\limits_{\mu:~ c_{\mu}>0}\P^{\mu}$ and
$\P^-:=\sum\limits_{\mu:~ c_{\mu}<0}\P^{\mu}$ be the
corresponding deformed symmetric and antisymmetric
projectors respectively. Assume that relations
(\ref{qcr1}), (\ref{qcr2}) define algebrae ${\cal M}$,
${\cal M}'$ with the same Poincar\'e series as their
classical counterparts. In order that   
relations (\ref{qcr1}-\ref{ansatz})
define a deformed Weyl algebra ${\cal A}_+^q$
(resp. Clifford algebra ${\cal A}_-^q$)
with the same Poincar\'e series as its classical counterpart 
${\cal A}_+$ (resp. ${\cal A}_-$) 
there must exist exactly one negative (resp. positive) $c_{\mu}$,
say $c_-$ (resp. $c_+$), and
the commutation relations (\ref{ansatz}) have to take one of the
two following forms
\eqa
A^i\, A^+_j &=& \delta^i_j{\bf 1} \:-\: (c_{\mp})^{-1}\R^{ih}_{jk}
A^+_h\, A^k, \label{first}\\
A^i\, A^+_j &=& \delta^i_j{\bf 1} \:-\: c_{\mp}\R^{-1}{}^{ih}_{jk}
A^+_h\, A^k. \label{second}
\ena
\label{lemma}
\end{lemma}
$Proof.$  Let us multiply eq. (\ref{qcr1}) by $A^l$ from the left. 
We easily find
\eqa
0 & = &A^l\P^{\mp~ij}_{~~hk}A^+_i\, A^+_j \cr
  & \stackrel{(\ref{ansatz})}{=} &
[\P^{\pm}({\bf
1}+S)]^{li}_{hk}A^+_i+(S_{12}S_{23}\P^{\pm}_{12})^{lij}_{hkm}
A^+_iA^+_jA^m \nonumber 
\ena
In order that the second term vanishes without introducing new,
third degree relations (which would yield a different Poincar\'e
series) it must be either $S\propto\R$ or 
$S\propto\R^{-1}$, so that
\[
(S_{12}S_{23}\P^{\pm}_{12})^{lij}_{hkm} A^+_iA^+_jA^m 
   \stackrel{(\ref{braid2})}{=} 
(\P^{\pm}_{23}S_{12}S_{23})^{lij}_{hkm} A^+_iA^+_jA^m 
   \stackrel{(\ref{qcr1})}{=} 0.
\]
These correspond to the two possible braidings \cite{majid}. 
If $S=b\R$, then the first term vanishes iff
\[
0=\P^{\pm}({\bf 1}+S)=\sum\limits_{\mu:~ \pm c_{\mu}>0}\P^{\mu}
({\bf 1}+c_{\mu}b)\qquad\Leftrightarrow \qquad {\bf 1}+c_{\mu}b
\x\x\forall\mu: \pm c_{\mu}>0.
\]
Thus  there must
exist exactly one positive (resp. negative)
$c_{\mu}$ and relation (\ref{first}) must hold.
Similarly one proves relation (\ref{second}) if
$S=b\R^{-1}$. $\Box$

As immediate consequences of this lemma and of the decompositions
(\ref{defproj}) we find:
\begin{itemize}
\item there exist no satisfactory definitions of the $q$-deformed
algebrae ${\cal A}^q_{-,SO(N),\phi_d}$, 
${\cal A}^q_{+,sp({n}),\phi_d}$, since these
correspond respectively
to the projectors $(\ref{defproj})_2$, $(\ref{defproj})_6$;
\item there exist satisfactory definitions of the $q$-deformed
algebrae 
 ${\cal A}^q_{+,SL(N),\phi_d}$ \cite{puwo,wezu}, ${\cal
A}^q_{-,SL(N),\phi_d}$ 
\cite{pusz},
${\cal A}^q_{+,SO(N),\phi_d}$ \cite{cawa}, ${\cal
A}^q_{-,sp({n}),\phi_d}$,
since these are the algebrae corresponding respectively
to the projector $(\ref{defproj})_4$  
$(\ref{defproj})_1$ , 
$(\ref{defproj})_5$ , $(\ref{defproj})_3$ 
(up to our knowledge, the latter
has never been considered before in the literature). 

\end{itemize}

\section{Main embedding prescription}
\label{qcrsec}

We would like to generalize the construction of the preceding section
to the case in which $A^+_i$, $A^i$ belong respectively to
corepresentations  
$\phi_M,\phi_M^\vee$ that are direct sums of $M\ge 1$ copies
of $\phi_d,\phi^\vee_d$.
Let $^{\alpha}{\cal A}^q_{\pm,G,\phi_d}$ 
($\alpha=1,2,...,M$) be $G_q$-covariant 
$q$-deformed Heisenberg algebra  with generators 
$\1, A^{\alpha,i},A^+_{\alpha,i}$, $i=1,....,N$,
and relations
\eqa
&&\P^{(\alpha)}{}_{ij}^{hk} A^+_{\alpha,h}A^+_{\alpha,k} =
 0\label{pippo1.1} \\
&&\P^{(\alpha)}{}^{ij}_{hk} A^{\alpha,k}A^{\alpha,h} = 0 
\label{pippo1.2}\\
&&A^{\alpha,i}A^+_{\alpha,j}-\delta^i_j\1-
(-1)^{\epsilon_{\alpha}}
[(q^{1-2\epsilon_{\alpha}}\R)^{\eta_{\alpha}}]^{ih}_{jk} 
A^+_{\alpha,h}A^{\alpha,k} = 0.\label{pippo1.3}
\ena
According to the last remark in the previous
section, let $\epsilon_{\alpha}$ take 
the values $\epsilon_{\alpha}\equiv 0$ if $G=SO(N)$, 
$\epsilon_{\alpha}\equiv 1$ if $G=Sp(n)$, and
$\epsilon_{\alpha}=0,1$ if $G=SL(N)$; $\epsilon_{\alpha}=0,1$ 
correspond to Weyl, Clifford respectively. Moreover, let
\eq 
\P^{(\alpha)}=\cases{\P^+\qquad 
\mbox{  if  }\epsilon_{\alpha}=0\cr 
\P^- \qquad \mbox{  if  }\epsilon_{\alpha}=1; \cr}
\en

Recalling that the  comodules of $\fun$ belong  to a 
braided monoidal
category, we know that consistent commutation 
relations between 
the generators of ~$^{\alpha}{\cal A}^q_{\pm,G,\phi_d}$,
$^{\beta}{\cal A}^q_{\pm,G,\phi_d}$,
$\alpha\neq \beta$, are
given by the two 
possible braidings (the latter correspond to the quasitriangular
structures ${\cal R},{\cal R}^{-1}_{21}$ \cite{majid}). 
Accordingly, the commutation relations between 
$A^+_{\alpha,i},A^+_{\beta,j}$ for instance may become
\eqa
\mbox{either} \qquad\qquad
A^+_{\alpha,i},A^+_{\beta,j} & \propto & \R^{hk}_{ij} 
A^+_{\beta,h}A^+_{\alpha,k}  \nonumber\\
\mbox{or} \qquad\qquad
A^+_{\alpha,i},A^+_{\beta,j} & \propto & \R^{-1}{}^{hk}_{ij} 
A^+_{\beta,h}A^+_{\alpha,k}.\nonumber
\label{cross1}
\ena
There are $\frac{M(M-1)}2$ different pairs $(\alpha,\beta)$;
if we could choose for each pair the upper or lower solution
independently
we would have $2^{\frac{M(M-1)}2}$ different versions of the deformed
commutation relations. We claim that, in fact, only $M!$ are allowed,
in other words that, up to a reordering (\ie a permutation of the
$\alpha$'s),
the only consistent way is:
\begin{prop}
Without loss of generality we can
assume
\eq
A^+_{\alpha,i},A^+_{\beta,j}=(-1)^{\epsilon_{\alpha}
\epsilon_{\beta}}
c_{\alpha\beta}\R^{hk}_{ij}
A^+_{\beta,h}A^+_{\alpha,k} ,  
\qquad\qquad \mbox{if~~} \alpha<\beta,\label{pippo2}
\en
with $c_{\alpha\beta}\stackrel{q\rightarrow 1}{\rightarrow}1$.
\end{prop}
(We have factorized the overall sign necessary 
to get the correct commutation relations between fermionic or
bosonic variables in the classical limit).

$Proof$. The claim can be proved inductively.
It is obvious if $M=2$. Assume now that the claim is true when
$M=P$, and call $A^+_{\cdot,i}$ the generators of
the $(P\!+\!1)$-th additional subalgebra. We need to prove that
\eqa
A^+_{\beta,i}\,A^+_{\cdot,j}\propto \R^{hk}_{ij}
A^+_{\cdot,h}\,A^+_{\beta,k} \qquad & \Rightarrow & \qquad
A^+_{\alpha,i}\,A^+_{\cdot,j}\propto \R^{hk}_{ij}
A^+_{\cdot,h}\,A^+_{\alpha,k} \x\x \forall\, \alpha<\beta \cr
A^+_{\gamma,i}\,A^+_{\cdot,j}\propto \R^{-1}{}^{hk}_{ij}
A^+_{\cdot,h}\,A^+_{\gamma,k} \qquad & \Rightarrow & 
\qquad
A^+_{\delta,i}\,A^+_{\cdot,j}\propto \R^{-1}{}^{hk}_{ij}
A^+_{\cdot,h}\,A^+_{\delta,k} \x\x \forall \,\delta>\gamma.
\nonumber
\ena
Let $A^+_{\beta,i}\,A^+_{\cdot,j}=V^{hk}_{ij}
A^+_{\cdot,h}\,A^+_{\beta,k}$;
we can invert the order of the factors in the
product $A^+_{\alpha,h}A^+_{\beta,i}\,A^+_{\cdot,j}$ 
either by permuting the first two factors, then the last two, 
finally the first two again, or by permuting the last two factors, 
then the first two, finally the last two again;
to get the same result we need that
$\R_{12}V_{23}\R_{12}=\R_{23}V_{12}\R_{23}$. This
 equation is satisfied
iff $V\propto \hat R$. Thus we have proved the first 
implication.
Similarly one proves the second. $\Box$

Eq. (\ref{pippo2}) and the condition that
$A^{\alpha,i}$ are the dual generators of $A^+_{\alpha,i}$ implies
(for $\alpha< \beta$)
\eq
A^{\alpha,j}A^{\beta,i} =(-1)^{\epsilon_{\alpha}
\epsilon_{\beta}}
c_{\alpha\beta}\R^{ij}_{hk}A^{\beta,k}A^{\alpha,h}
\label{pippo3.2}
\en
As for the remaining relations, we shall look for them in the
form 
$A^{\beta,i}A^+_{\alpha,j} = M^{ih}_{jk} A^+_{\alpha,h}A^{\beta,k}$.
It is easy to check that from either of the previous relation and
the commutation relations of $^{\alpha}{\cal A}^q_{\pm,G,\phi_d}$
it follows  (for $\alpha< \beta$):
\eqa
A^{\beta,i}A^+_{\alpha,j}&=&(-1)^{\epsilon_{\alpha}
\epsilon_{\beta}}
c_{\alpha\beta}\R^{ih}_{jk} A^+_{\alpha,h}A^{\beta,k}\label{pippo3.1}\\
A^{\alpha,i}A^+_{\beta,j}&=&(-1)^{\epsilon_{\alpha}
\epsilon_{\beta}}
c_{\alpha\beta}^{-1}(\R^{-1})^{ih}_{jk} A^+_{\beta,h}
A^{\alpha,k}\label{pippo3.3}
\ena
For instance,
relation (\ref{pippo3.1}) is derived  by consistency when requiring 
that one gets the same result from
$A^{\alpha,i}A^+_{\alpha,j}A^+_{\beta,k}$
either by permuting the first two factors, then the last two, 
finally the first two again, or by permuting the last two factors, 
then the first two, finally the last two again.

We will call ${\cal A}^q_{G,\phi_M}$
the unital algebra generated by 
$\1,A^{\alpha,i},A^+_{\alpha,i}$, 
$\alpha=1,2,...,M$, $i=1,...,N$ and commutation relations
(\ref{pippo1.1}-\ref{pippo1.3}), 
(\ref{pippo2}-\ref{pippo3.3}). 
We have thus proved 

\begin{prop}
 ${\cal A}^q_{G,\phi_M}$ has the same Poincar\'e series as its classical
counterpart  ${\cal A}_{G,\phi_M}$.
\end{prop}

\renewcommand{\theequation}{\thesection.\arabic{equation}}
\sect{$*$-Structures}
\label{star}

Let $\fun$ be a Hopf $*$-algebra  
and assume that 
${}^{\alpha}{\cal A}^q_{\pm,\footnotesize{G},\phi_d}$
are \fun-comodule $*$-algebrae:
\eq
\phi_d(b^{\star_{\alpha}})=[\phi_d(b)]^{\star_{\alpha}\otimes *}, 
\qquad\qquad \x\x b\in 
{}^{\alpha}{\cal A}^q_{\pm,{\footnotesize\mbox{\bf g\,}},\phi_d},
\label{condstar}
\en
(here ``$\star_{\alpha}$'' denotes the $*$ of 
${}^{\alpha}{\cal A}^q_{\pm,\footnotesize{G},\phi_d}$).
Can we use the $\star_{\alpha}$'s to 
build a $*$-structure $\star$ of the whole 
\Aq?

In the case that $*$ realizes the compact real section
of \fun (what requires $q\in\rn^+$), 
then the simplest $*$-structure in 
${\cal A}^q_{\pm,{\footnotesize\mbox{\bf g\,}},\phi_d}$
is
\eq
(A^i)^{\star}= A^+_i.
\label{star2}
\en

It is immediate to check that the Ansatz
$(A^{i,\alpha})^{\star}= A^+_{i,\alpha}$ would be
compatible with relations (\ref{pippo1.1}-\ref{pippo1.3}), but
inconsistent
with relations (\ref{pippo2}-\ref{pippo3.3}). 
Therefore let us choose the Ansatz
\eq
(A^{i,\alpha})^{\star}= A^+_{i,\pi(\alpha)},
\en
where $\pi$ is some permutation of $(1,\ldots,M)$.
It is easy to check that
consistency with relations (\ref{pippo2}-\ref{pippo3.3})
requires
\eqa
&& \pi(\alpha) = M-\alpha+1, \label{consi1}\\
&& \eta_{\pi(\alpha)}=\eta_{\alpha}, \qquad
c_{\pi(\alpha)\pi(\beta)}=c_{\beta\alpha},\qquad
\epsilon_{\pi(\alpha)}=\epsilon_{\alpha}. \label{consi2}
\ena
Eq. (\ref{consi1}) shows that $\pi$ must 
be the inverse-ordering permutation;
Eq. (\ref{consi2})$_3$ amounts to say that
$\star$ must preserve the bosonic or fermionic
character of the generators.

${\cal A}_{+,SO(N),\phi_d}^q$ admits also an alternative
$*$-structure compatible with $\phi_d$, namely
\eq
(A^+_i)^{\star}= A^+_jC_{ji},
\en
together with a nonlinear transformation for $(A^i)^{\star}$
\cite{oleg2}.
Here $C_{ij}$ is the $q$-deformed metric matrix \cite{frt},
which is related to the projector ${\cal P}^t$ appearing in
(\ref{decompo})$_2$ through the formula
${\cal P}^t{}^{ij}_{hk}=\frac{C^{ij}C_{hk}}{C^{lm}C_{lm}}$.
It is easy to check that the
Ansatz
\eq
(A^+_{i,\alpha})^{\star}= A^+_{j,\pi(\alpha)}C^{ji},
\en
together the corresponding 
nonlinear one for $(A^{i,\alpha})^{\star}$, defines a 
consistent $*$-structure of 
${\cal A}_{+,SO(N),\phi_d}^q$ provided 
relations (\ref{consi1}), (\ref{consi2}) hold (with 
$\epsilon_{\alpha}\equiv 0$ $\forall \alpha$).

\renewcommand{\theequation}{\thesection.\arabic{equation}}
\sect{Modified prescription: $GL_q(M)\times G_q$-covariant algebrae}
\label{glcov}

If all the generators of $\A$
have the same Grassman parity, they belong to a
corepresentation
of $GL(M)\times G$. The coaction of the group $GL(M)$
amounts
to a linear invertible transformation $T$ of the $a^{\alpha,i}$
and of the $a^+_{\alpha,i}$:
\eq
a^{\alpha,i}\rightarrow a^{\beta,i}T^{\alpha}_{\beta}\qquad\qquad
a^+_{\alpha,i}\rightarrow a^+_{\beta,i}T^{-1}{}_{\alpha}^{\beta},
\en
which leaves the commutation relations
(\ref{ccr1}-\ref{ccr3}) invariant.
[If in addition 
we require some $*$-structure to be preserved, then  $T$ has to belong
to
some suitable subgroup of $GL(M)$; \eg
$T\in U(M)$ if $(a^i)^{\dagger}=a^+_i$.]
We try to construct now a variant of the algebra of section
\ref{qcrsec}  
having explicitly $GL_q(M)\times G_q$-covariant 
generators\footnote{Or equivalently $SL_q(M)\times G_q$-covariance,
if we impose also the unit condition on the $q$-determinant of
$GL_q(M)$.}.

Let $T^{\alpha}_{\beta}$, 
$t=\mbox{det}_q\Vert T^{\alpha}_{\beta}\Vert$ be the generators of the 
quasitriangular Hopf algebra $Fun[GL_q(M)]$,  and $T^a_b$
the generators of $Fun(G_q)$ \cite{frt}.
Let us introduce collective indices $A,B,...$, denoting the pairs
$(\alpha,a)$, $(\beta,b),...$. 
The Hopf algebra $Fun(GL_q(M)\times G_q)$
can be defined as the algebra generated by objects $T^A_B$
satisfying commutation relations
which can be obtained from (\ref{rel1}) by the
replacement 
\eq
T^A_B\rightarrow T^{\alpha}_{\beta}T^a_b
\label{repla}
\en
by assuming that $[T^{\alpha}_{\beta},T^a_b]=0$:
\eq
\RH^{AB}_{CD}T^C_ET^D_F=T^A_CT^B_D\RH^{CD}_{EF}.
\en
Here $\RH$ is one of the matrices
\eq
\RH_{\pm}{}^{AB}_{CD}:=\R^{\pm 1}_M{}^{\alpha\beta}_{\gamma\delta}
\R^{ab}_{cd}\equiv (\R^{\pm 1}_M\ot
\R)^{AB}_{CD},
\en
and $\R_M$ is the braid matrix (\ref{explicit}) of $SL_q(M)$.
$\RH_{\pm}$ satisfies the braid equation, since $\R,\R_M$ do. 
The coproduct, counit, antipode and quasitriangular structure 
are introduced as in Sect. \ref{preli} by
$\Delta(T^A_B)=T^A_C\otimes T^C_B$, $\varepsilon(T^A_B)=\delta^A_B$,
$ST^A_B=T^{-1}{}^A_B$.

A (right) comodule algebra 
of $Fun(GL_q(M)\times G_q)$ can be associated to the defining 
corepresentation of the latter,
$\phi_D(A^+_A)=A^+_B\otimes T^B_A$, where $A^+_C$ denote the
generators. The dual comodule algebra, with generators $A^C$,
will be associated to the contragradient corepresentation
$\phi_D^\vee(A^A)=A^B\otimes ST^A_B$.
To find compatible quadratic commutation relations 
among the $A^+_B$'s (resp. $A^B$'s) we need
the projector decomposition of $\RH_{\pm}$, as in Sect. \ref{preli}.
For this scope we just need
to write down the projector decompositions of both $\R_M^{\pm 1}$ and 
$\R$ and note that 
$\PH:=\P_M\ot\P'$
is a projector $\PH$ whenever $\P,\P'$ are.

We start with the case $G_q=SL_q(N)$. We find
\eqa
\RH_+ & = & (q\P_M^S-q^{-1}\P_M^A)\ot (q\P^S-q^{-1}\P^A) \cr
      & = & -(\P_M^S\ot \P^A+\P^A_M\ot\P^S)+ q^2\P^S\ot\P^S
            +q^{-2} \P_M^A\ot \P^A_M \cr
      & =:& -\PH^-+q^2\PH^{S,1}+q^{-2}\PH^{S,2} 
\label{deco+}
\ena
and
\eqa
\RH_- & = &(q^{-1}\P^S_M-q\P_M^A)\ot (q\P^S-q^{-1}\P^A) \cr
      & = & (\P^A_M\ot \P^A+\P^S_M\ot\P^S)-q^2\P^A_M\ot\P^S
            -q^{-2} \P^S_M\ot \P^A \cr
      & =:&\PH^+-q^2\PH^{A,1}-q^{-2}\PH^{A,2}.
\label{deco-}
\ena

We are in the condition to apply Lemma \ref{lemma}.
As a consequence, there exists
a $GL_q(M)\times SL_q(N)$-covariant Weyl algebra
${\cal A}^q_{+,\footnotesize{GL_q(M)\times SL_q(N)},\phi_D}$, defined
by the following commutation relations:
\eqa
\PH^-{}^{CD}_{AB}A^+_CA^+_D & = & 0 \label{lilla1}\\
\PH^-{}_{CD}^{AB}A^DA^C & = & 0 \label{lilla2}\\
A^AA^+_B-\delta^A_B\1-\RH_+{}^{AC}_{BD}A^+_C A^D & = & 0.
\label{lilla3}
\ena
Moreover, there exists
a $q$-deformed
$SL_q(M)\times SL_q(N)$-covariant Clifford algebra
${\cal A}^q_{-,\footnotesize{SL_q(M)\times SL_q(N)},\phi_D}$, defined
by the following commutation relations:
\eqa
\PH^+{}^{CD}_{AB}A^+_CA^+_D & = & 0\label{lilla4}\\
\PH^+{}_{CD}^{AB}A^DA^C & = & 0 \label{lilla5}\\
A^AA^+_B-\delta^A_B\1+\RH_-{}^{AC}_{BD}A^+_C A^D & = & 0
\label{lilla6}
\ena
According to lemma \ref{lemma}, one could give also 
alternative definitions 
with $\RH^{-1}$ instead of $\RH$ in relations (\ref{lilla3}),
(\ref{lilla6}).

Let us verify that relations (\ref{lilla1}), (\ref{lilla2}),
(\ref{lilla4}), (\ref{lilla5}) are of the kind considered in section
\ref{qcrsec}.

We take first relations (\ref{lilla1}) into account. We find
\eqa
(q+q^{-1})^2\PH^- & \stackrel{(\ref{deco+})}{=} & (q+q^{-1})^2
[\P^S\ot \P^A+\P^A+\ot \P^S] \nonumber\cr
 & \stackrel{(\ref{defproj})_4}{=} &
(q\1- \R_M)\ot (q^{-1}\1+\R)+ (q^{-1}\1+\R_M) \ot (q\1- \R) \nonumber\cr
 & = & 2(\1\ot \1- \R_M\ot \R)+ (q-q^{-1})(\1\ot \R+ \R_M\ot\1).
\ena
Using relation (\ref{explicit}) we can write $\R_M$ explicitly
and check that relations  (\ref{lilla1})
amount to relations 
\eqa
\P^-{}_{ij}^{hk} A^+_{\alpha,h}A^+_{\alpha,k} &=
 &0, \label{zaz1} \\
A^+_{\alpha,i},A^+_{\beta,j} -\R^{hk}_{ij}
A^+_{\beta,h}A^+_{\alpha,k} & = & 0,  
\qquad\qquad \mbox{if~~}\alpha<\beta \label{zaz2}.
\ena
Similarly one verifies that: 1) relations (\ref{lilla2}) amount
to relations 
\eqa
\P^-{}^{ij}_{hk} A^{\alpha,k}A^{\alpha,h} &= &0 \label{zaz3}\\
A^{\alpha,j}A^{\beta,i}-\R^{ij}_{hk}A^{\beta,k}A^{\alpha,h}& = & 0
\qquad\qquad \mbox{if~~}\alpha<\beta; \label{zaz4}
\ena
2) that relations (\ref{lilla4}) amount to 
relations\footnote{These are of the type considered in Sect. 2,
provided we invert the order of greek indices.}
\eqa
\P^+{}_{ij}^{hk} A^+_{\alpha,h}A^+_{\alpha,k} &= &0\label{zaz5} \\
A^+_{\alpha,i},A^+_{\beta,j} +\R^{-1}{}^{hk}_{ij}
A^+_{\beta,h}A^+_{\alpha,k} & = & 0 ,  
\qquad\qquad \mbox{if~~}\alpha<\beta;\label{zaz6}
\ena
3) that relations (\ref{lilla5}) amount to relations 
\eqa
\P^+{}^{ij}_{hk} A^{\alpha,k}A^{\alpha,h} &= &0 \label{zaz7}\\
A^{\alpha,j}A^{\beta,i} + \R^{-1}{}^{ij}_{hk}A^{\beta,k}A^{\alpha,h}
& = & 0\qquad\qquad \mbox{if~~}\alpha<\beta .\label{zaz8}
\ena

On the other hand,
relations (\ref{lilla3}), (\ref{lilla6}) for $\alpha\neq \beta$
are {\it not} of the type (\ref{pippo3.1}), (\ref{pippo3.3})
found in section \ref{qcrsec}; in fact, in
a similar way one can show that relation (\ref{lilla3})
takes the form
\eqa
&&A^{\alpha,a}A^+_{\beta,b} - \R^{ac}_{bd}A^+_{\beta,c}A^{\alpha,d} = 0
\z\z\z\z\z \alpha\neq\beta, \z\label{d+}\\
&&A^{\alpha,a}A^+_{\alpha,b} - \delta^a_b{\bf 1} -
q\R^{ac}_{bd}A^+_{\alpha,c}A^{\alpha,d} - (q-q^{-1})
\sum\limits_{\beta>\alpha}\R^{ac}_{bd}A^+_{\beta,c}A^{\beta,d} = 0\z
\label{s+}
\ena
whereas relation (\ref{lilla6}) amounts to
\eqa
&&A^{\alpha,a}A^+_{\beta,b} + \R^{ac}_{bd}A^+_{\beta,c}A^{\alpha,d} = 0
\z\z\z\z\z\alpha\neq \beta, \z \label{d-}\\
&&A^{\alpha,a}A^+_{\alpha,b} - \delta^a_b {\bf 1} + 
q^{-1}\R^{ac}_{bd}A^+_{\alpha,c}A^{\alpha,d}- (q-q^{-1})
\sum\limits_{\beta<\alpha}\R^{ac}_{bd}A^+_{\beta,c}A^{\beta,d} = 0\z.
\label{s-}
\ena
Relations (\ref{d+}), (\ref{d-}) specialized to the case $\alpha>\beta$
coincide with relations (\ref{pippo3.1}); specialized to the case
$\alpha<\beta$, they differ from relations (\ref{pippo3.3}). 
Relations (\ref{s+}), (\ref{s-})
differ from relations (\ref{pippo1.3}) by the additional terms 
with coefficient $(q-q^{-1})$.

The subalgebra ${\cal M}$
(resp. ${\cal M}'$) generated by $A^+_A$'s (resp.
$A^A$'s) has the same Poincar\'e series of the 
subalgebra generated by classical $a^+_{\alpha a}$'s (resp. $a^{\alpha
a}$'s),
because of relations (\ref{zaz1}), (\ref{zaz2}) 
[resp. (\ref{zaz3}), (\ref{zaz4})] in the Weyl case and
because of relations (\ref{zaz5}), (\ref{zaz6}) 
[resp. (\ref{zaz7}), (\ref{zaz8})] in the Clifford case.
Since relations (\ref{lilla3}), (\ref{lilla6}) allow to change the order
of $A^+_A$'s and $A^B$'s in any product, we conclude that
\begin{prop}
The algebrae
${\cal A}^q_{\pm,\footnotesize{GL_q(M)\times SL_q(N)},\phi_D}$ 
have the same Poincar\'e series as their classical counterparts.
\end{prop}

Finally, let us ask about $*$-structures.
When $q\in\rn^+$ the Hopf algebra $GL_q(M)\times SL_q(N)$
admits the compact section $U_q(M)\times SU_q(N)$ \cite{frt}.
The deformed Heisenberg
algebrae defined by relations (\ref{lilla1}-\ref{lilla6})
admit a natural $U_q(M)\times SU_q(N)$-covariant
$*$-structure given by
\eq
(A^A)^{\star}=A^+_A;
\en
this can be easily checked by applying this $\star$  to relations 
(\ref{lilla1}-\ref{lilla6}) and by
noting that $\R^T=\R$ and therefore
$\RH^T=\RH$, $\PH^T=\PH$.

\medskip

Let us take now in consideration
the cases that $G_q=SO_q(N),Sp_q(n)$. The projector
decomposition of $\R_M\ot \R=\sum_{\mu} \lambda_{\mu}\PH^{\mu}$
gives $\lambda_{\mu}=q^2, q^{-2}, -1,\pm q^{2-N}, \mp q^{-N}$,
where the
upper and lower sign refer to $G_q=SO_q(N)$ and $Sp_q(n)$ respectively.
  The projector
decomposition of $\R^{-1}_M\ot \R=\sum_{\mu} \lambda_{\mu}\PH^{\mu}$
gives $\lambda_{\mu}=-q^2, -q^{-2}, 1,\mp q^{2-N}$, $\pm q^{-N}$.
In both cases we always have more
than one positive and more than one negative $\lambda_{\mu}$.  
By Lemma \ref{lemma} no
$GL_q(M)\times G_q$ covariant $q$-deformed Weyl/Clifford algebra can
be built by this procedure.

\section*{Acknowledgments}

It is a pleasure to thank J.\ Wess  for his scientific
support and for the 
hospitality at his Institute.
This work was financially supported through a TMR fellowship
granted by the European Commission, Dir. Gen. XII for Science,
Research and Development, under the contract ERBFMICT960921.


\end{document}